\font \bigbf=cmbx10 scaled \magstep1

\font\math=msbm10
 
\def\paragraphe#1{\bigskip\goodbreak {\bigbf #1}\par\nobreak\vskip
12pt plus 2pt\nobreak} 
\def\alinea#1{\medskip\allowbreak{\bf#1}\par\nobreak\vskip 9pt plus
1pt\nobreak} 
\def\th#1{\bigskip\goodbreak {\bf Theor\`eme #1.} \par\nobreak \sl }
\def\prop#1{\bigskip\goodbreak {\bf Proposition #1.} \par\nobreak \sl }
\def\lemme#1{\bigskip\goodbreak {\bf Lemme #1.} \par\nobreak \sl }
\def\cor#1{\bigskip\goodbreak {\bf Corollaire #1.} \par\nobreak \sl }
\def\dem{\bigskip\goodbreak \it D\'emonstration. \rm}
\def\ndem{\bigskip\goodbreak \rm}
\def\qed{\par\nobreak\hfill $\bullet$ \par\goodbreak}
\def\uple#1#2{#1_1,\ldots ,{#1}_{#2}}
\def\corde#1#2-#3{{#1}_{#2},\ldots ,{#1}_{#3}}

\def\ordcorde#1#2-#3{{#1}_{#2} \le \cdots \le {#1}_{#3}}
\def\strictordcorde#1#2-#3{{#1}_{#2} < \cdots < {#1}_{#3}}
\def \restr#1{\mathstrut_{\bigl|}\raise-8pt\hbox{$#1$}}
\def \inver{^{-1}}
\def\fleche#1{\mathop{\hbox to #1 mm{\rightarrowfill}}\limits}
\def\gfleche#1{\mathop{\hbox to #1 mm{\leftarrowfill}}\limits}
\def\inj#1{\mathop{\hbox to #1
mm{$\lhook\joinrel$\rightarrowfill}}\limits} 
\def\ginj#1{\mathop{\hbox to #1
mm{\leftarrowfill$\joinrel\rhook$}}\limits} 
\def\surj#1{\mathop{\hbox to #1 mm{\rightarrowfill\hskip
2pt\llap{$\rightarrow$}}}\limits} 
\def\gsurj#1{\mathop{\hbox to #1 mm{\rlap{$\leftarrow$}\hskip 2pt
\leftarrowfill}}\limits}

\def \mop#1{\mathop{\hbox{\rm #1}}\nolimits}
\def \smop#1{\mathop{\hbox{\sevenrm #1}}\nolimits}

\def \bib #1{\null\medskip\hskip -14mm\hbox to 14mm{[#1]\hfill}}
\def\A{{\cal A}}

\def\P{{\cal P}}
\def\J{{\cal J}}

\def\H{{\cal H}}
\magnification=1200
\parindent=0cm
\centerline{\bigbf L'ALGEBRE DE HOPF BITENSORIELLE}
\vskip 12mm plus 2mm
\centerline{\bf Dominique Manchon}
\bigskip
\centerline{Institut Elie Cartan}
\centerline{CNRS}
\centerline{BP 239}
\centerline{54506 Vandoeuvre les Nancy Cedex, France}
\centerline{e-mail~: manchon@iecn.u-nancy.fr}
\vskip 12mm plus 2mm
\centerline{\bf R\'esum\'e}
\medskip\qquad 
Nous construisons pour tout corps $k$ de caract\'eristique z\'ero un
foncteur de la cat\'egorie des $k$-espaces vectoriels dans la
cat\'egorie des $k$-alg\`ebres de Hopf point\'ees, qui \`a tout espace
vectoriel $V$ associe son {\sl alg\`ebre de Hopf bitensorielle
point\'ee\/} $\A_V$. Cette alg\`ebre de Hopf est gradu\'ee, v\'erifie
une propri\'et\'e universelle, et contient une famille remarquable
d'\'el\'ements primitifs $\P$. Nous conjecturons que $\P$ engendre
l'alg\`ebre de Lie des \'el\'ements primitifs de $\A_V$. Enfin lorsque
$V$ est de dimension finie nous mettons en \'evidence un couplage de
Hopf entre $\A_V$ et $\A_{V^*}$ dont le noyau contient l'id\'eal (de
Hopf) engendr\'e par les \'el\'ements de $\P$ de degr\'e au moins
\'egal \`a $2$.  
\bigskip
\centerline{\bf Abstract}
\medskip\qquad
For any field $k$ of zero characteristic we give a functor from the
category of $k$-vector spaces into the category of $k$-Hopf algebras,
attaching to any vector space $V$ its {\sl bitensorial pointed Hopf
algebra\/} $\A_V$. This Hopf algebra is graded, fulfills a universal
property, and contains a remarkable subspace $\P$ of primitive
elements, which as a conjecture may generate the Lie algebra
$\mop{Prim}\A_V$. In case $V$ is finite-dimensional we exhibit a Hopf
pairing between $\A_V$ and $\A_{V^*}$ whose kernel contains the (Hopf)
ideal generated by the elements of $\P$ of degree $\ge 2$. 
\paragraphe{Introduction}
\qquad
Etant donn\'e un espace vectoriel $V$ sur un corps $k$, on sait lui
associer de mani\`ere fonctorielle une alg\`ebre, son alg\`ebre
tensorielle $T(V)$, qui contient $V$ de mani\`ere naturelle, et qui
v\'erifie de plus la propri\'et\'e universelle suivante : pour toute
alg\`ebre $A$ et pour toute application lin\'eaire $f$ de $V$ dans $A$
il existe un unique morphisme d'alg\`ebres $\overline f$ de $T(V)$
dans $A$ qui prolonge $f$. 
\smallskip\qquad
Nous montrons ici qu'une construction similaire est possible dans la
cat\'egorie des alg\`ebres de Hopf point\'ees~: on consid\`ere sur
l'espace $T(V)$ la structure de cog\`ebre obtenue en dualisant la
structure d'alg\`ebre de $T(V^*)$. La comultiplication est donn\'ee
par~: 
$$\Delta(x_1\otimes\cdots\otimes x_n)=\sum_{i=0}^n
(x_1\otimes\cdots\otimes x_j)\tilde\otimes (x_{j+1}\otimes \cdots
\otimes x_n)$$ 
et la co-unit\'e $\varepsilon$ par le terme constant. Appliquant le
foncteur "alg\`ebre tensorielle" \`a la comultiplication $\Delta$
et \`a la co-unit\'e $\varepsilon$ nous
munissons de mani\`ere naturelle la double alg\`ebre tensorielle
$T(T(V))$ d'une structure de big\`ebre. Identifiant les puissances de
l'unit\'e $c$ de $T(V)$ et l'unit\'e $1$ de $T(T(V))$, qui sont les
\'el\'ements de type groupe de la big\`ebre, on obtient une big\`ebre
point\'ee $\A_V$ qui se trouve \^etre une alg\`ebre de Hopf. 
\smallskip\qquad
Lorsque le corps est de caract\'eristique nulle, l'antipode est
d'ordre infini et admet une expression explicite~: on d\'efinit $S_0$
comme l'unique antimorphisme d'alg\`ebres tel que~: 
$$S_0\restr{T(V)}=-I+2u\varepsilon$$
et l'{\sl op\'erateur de c\'esure\/} $U$ comme l'unique d\'erivation
de $\A_V$ telle que~: 
$$U\restr{T(V)}=(I-u\varepsilon)*(I-u\varepsilon)$$
o\`u l'\'etoile d\'esigne le produit de convolution [Ab, Sw]. On a
alors~: (Th\'eor\`eme I.3)~: 
$$S=(\exp -U).S_0$$
\qquad La construction est fonctorielle, et l'alg\`ebre de Hopf
point\'ee $\A_V$ ainsi construite v\'erifie la propri\'et\'e
universelle suivante~: pour tout espace vectoriel $V$ sur $k$, pour
toute alg\`ebre de Hopf point\'ee $\H$ et pour tout morphisme de
cog\`ebres $f$ de $T(V)$ dans $\H$ il existe un unique morphisme
d'alg\`ebres de Hopf $\overline f$ de $\A_V$ dans $\H$ qui prolonge
$f$. On peut enlever le mot "point\'ee" dans l'\'enonc\'e en rempla\c
cant $\A_V$ par $\tilde \A_V$, obtenue \`a partir de $T(T(V))$ en
rajoutant formellement les inverses des \'el\'ements de type groupe
(\S\ I.6). 
\smallskip\qquad
Le deuxi\`eme r\'esultat est le th\'eor\`eme I.5, qui met en
\'evidence une famille d'\'el\'ements primitifs~: on montre que pour
tout tenseur sym\'etrique $v\in T(V)$ l'\'el\'ement $\varphi(U)v$ est
primitif, o\`u $\varphi$ est la s\'erie enti\`ere d\'efinie par~: 
$$\varphi(z)={1-e^{-z}\over z}$$
\qquad
Nous conjecturons que cette famille engendre l'alg\`ebre de Lie des
\'el\'ements primitifs de $\A_V$. 
\smallskip\qquad Dans la deuxi\`eme partie nous mettons en \'evidence,
dans le cas o\`u l'espace vectoriel $V$ est de dimension finie, un
couplage de Hopf entre $\A_V$ et $\A_{V^*}$, dont le noyau contient
l'id\'eal engendr\'e par l'ensemble des $\varphi(U)v,\, v\in
S^{(2)}(V)$. (Th\'eor\`eme II.1). 
\smallskip\qquad
On obtient bon nombre d'alg\`ebres int\'eressantes (alg\`ebre
sym\'etrique ou ext\'erieure, alg\`ebre enveloppante d'une alg\`ebre
de Lie, etc.) comme quotient de l'alg\`ebre tensorielle par un id\'eal
bien choisi. On peut de m\^eme esp\'erer retrouver des alg\`ebres de
Hopf comme quotients ou sous-quotients de l'alg\`ebre de Hopf
bitensorielle par un id\'eal de Hopf ad hoc. Par exemple les
alg\`ebres enveloppantes quantifi\'ees pourraient \^etre obtenues \`a
partir de l'alg\`ebre de Hopf topologique $\A_V[[h,t]]$ sur l'anneau
$k[[h,t]]$, o\`u $V$ est une big\`ebre de Lie [Dr,E-K, R]. Nous
esp\'erons revenir sur cette question dans un prochain article. 
\medskip\qquad
L'auteur tient \`a remercier Alain Fuser, Pierre-Yves Gaillard, Guy
Rousseau et Marc Rosso pour d'utiles discussions, ainsi que Pierre
Marchand pour les aspects combinatoires du th\'eor\`eme I.5, et
Olivier Ramar\'e pour la d\'emonstration du non moins combinatoire
lemme II.7.  
\paragraphe{I. L'alg\`ebre de Hopf bitensorielle}
\qquad Soit $V$ un espace vectoriel de dimension finie sur un corps
$k$, et soit $T(V)$ son alg\`ebre tensorielle. On munit un quotient de
la double alg\`ebre tensorielle $T(T(V))$ d'une structure naturelle
d'alg\`ebre de Hopf, ni commutative ni cocommutative, dont l'antipode,
d'ordre infini, admet une expression explicite. 
\alinea{I.1. La big\`ebre bitensorielle}
On d\'esigne par $T(V)$ l'alg\`ebre tensorielle de $V$ d\'efinie par~:
$$T(V)=kc\oplus V\oplus V^{\otimes 2}\oplus\cdots$$
o\`u $c$ d\'esigne l'unit\'e du corps $k$. C'est l'\'el\'ement central
normalis\'e de l'alg\`ebre $T(V)$. On s'int\'eressera \`a la structure
de cog\`ebre sur $T(V)$ dont la comultiplication~: 
$$\Delta:T(V)\longrightarrow T(V)\tilde\otimes T(V)$$
est d\'efinie par $\Delta(c)=c\tilde\otimes c$ et~:
$$\Delta(x_1\otimes\cdots\otimes x_k)=\sum_{j=0}^k
(x_1\otimes\cdots\otimes x_j)\tilde\otimes(x_{j+1}\otimes\cdots\otimes
x_k)$$ 
On a not\'e $\tilde\otimes$ le produit tensoriel ci-dessus pour ne pas
le confondre avec le produit tensoriel $\otimes$ constitutif de
$T(V)$. On v\'erifie facilement que $\Delta$ est co-associative, et
que l'application de $T(V)$ dans $k$ qui \`a tout \'el\'ement associe
son "terme constant" est la co-unit\'e. 
\smallskip\qquad
On consid\`ere alors la double alg\`ebre tensorielle~:
$${\cal A}_0=T(T(V))=k\oplus T(V)\oplus T(V)^{\bullet 2}\oplus\cdots$$
o\`u l'on note avec un gros point ($\bullet$) le second produit
tensoriel, pour ne pas le confondre avec le $\otimes$ de $T(V)$. La
comultiplication se prolonge en un unique morphisme d'alg\`ebres~: 
$$\Delta:{\cal A}_0\longrightarrow T(T(V)\tilde\otimes T(V))$$
Or l'alg\`ebre $T(T(V)\tilde\otimes T(V))$ se plonge dans ${\cal
A}_0\tilde\otimes{\cal A}_0$ via~: 
$$(a_1\tilde\otimes b_1)\bullet\cdots\bullet (a_k\tilde\otimes
b_k)\longmapsto 
(a_1\bullet\cdots\bullet a_k)\tilde\otimes (b_1\bullet\cdots\bullet
b_k)$$ 
De m\^eme la co-unit\'e s'\'etend en un unique morphisme d'alg\`ebres~:
$$\varepsilon:{\cal A}_0\longrightarrow k$$
ce qui munit $({\cal A}_0,\bullet,\Delta,1,\varepsilon)$ d'une
structure naturelle de big\`ebre. 
\smallskip\qquad On appelle {\sl mot} un \'el\'ement ind\'ecomposable
de $T(V)$, et on appelle {\sl phrase} un \'el\'ement ind\'ecomposable
$v_1\bullet\cdots\bullet v_k$ de $\A_0$, o\`u chaque $v_j$ est
lui-m\^eme un mot. 
\alinea{I.2. L'alg\`ebre de Hopf bitensorielle point\'ee}
On consid\`ere dans ${\cal A}_0$ l'id\'eal bilat\`ere ${\cal I}$
engendr\'e par $c-1$. La relation~: 
$$\Delta(c-1)=c\tilde\otimes c-1\tilde\otimes 1=(c-1)\tilde\otimes
c+1\tilde\otimes(c-1)$$ 
montre que ${\cal I}$ est un bi-id\'eal. Le quotient~:
$${\cal A}={\cal A}_0/{\cal I}$$
h\'erite donc naturellement d'une structure de big\`ebre [Ab
th. 4.2.1]. Les \'el\'ements de type groupe de $\A_0$ (c'est \`a dire
v\'erifiant $\Delta(x)=x\tilde\otimes x$) sont les puissances de
$c$. Le seul \'el\'ement de type groupe de $\A$ est donc l'unit\'e,
c'est \`a dire que la big\`ebre $\A$ est point\'ee. Par ailleurs le
plongement canonique de $T(V)$ dans son alg\`ebre tensorielle $\A_0$
induit un plongement naturel de $T(V)$ dans $\A$, et $T(V)$ engendre
alors l'alg\`ebre $\A$. 
\smallskip\qquad
On note $m$ la multiplication de ${\cal A}\tilde\otimes{\cal A}$ dans
${\cal A}$, et $u$ l'application unit\'e~: $k\rightarrow {\cal A}$ qui
\`a $\lambda$ associe $\lambda.1$. On cherche un antipode, c'est \`a
dire par d\'efinition une application $S:{\cal A}\rightarrow {\cal A}$
telle que le diagramme ci-dessous commute~: 
$$(*)\hbox to 1cm{}\matrix{&&\A\tilde\otimes \A &
\fleche{8}^{\textstyle S\tilde\otimes I} &\A\tilde\otimes \A &&\cr 
	&\mathop{\nearrow}\limits^{\textstyle
	\Delta}&&&&\mathop{\searrow}\limits^{\textstyle m}&\cr 
\A&&\fleche {12}^{\textstyle \varepsilon}&k&\fleche {12}^{\textstyle
u} &&\A\cr &\mathop{\searrow}\limits_{\textstyle\Delta}
&&&&\mathop{\nearrow}\limits_{\textstyle m}&\cr 
&&\A\tilde\otimes \A &\fleche{8}_{\textstyle I\tilde\otimes S}
&\A\tilde\otimes \A &&\cr}$$ 
On a $\Delta(1)=1\tilde\otimes 1$ d'o\`u n\'ecessairement
$S(1)=1$. Pour tout $x\in V$ on a~: $\Delta(x)=1\tilde\otimes
x+x\tilde\otimes 1$ (puisque $c$ et $1$ sont maintenant identiques
dans $\A$), d'o\`u on d\'eduit, compte tenu du diagramme ci-dessus~: 
$$S(x)=-x$$
Essayons de d\'eterminer $S(x\otimes y)$ pour $x,y\in V$~: la partie
sup\'erieure du diagramme se lit~: 
$$m(I\tilde\otimes S)\Delta(x\otimes y)=0$$
soit~:
$$m(I\tilde\otimes S)(1\tilde\otimes\ x\otimes y+x\tilde\otimes
y+x\otimes y\ \tilde\otimes 1)=0$$ 
soit encore~:
$$S(x\otimes y)-x\bullet y+x\otimes y=0$$
d'o\`u finalement~:
$$S(x\otimes y)=x\bullet y-x\otimes y$$
On v\'erifie imm\'ediatement que la partie inf\'erieure du diagramme
commute, c'est \`a dire que~: 
$$m(S\tilde\otimes I)\Delta(x\otimes y)=0$$
Au cran suivant, en suivant la m\^eme m\'ethode on trouve~:
$$S(x\otimes y\otimes z)=-x\otimes y\otimes z+x\bullet(y\otimes
z)+(x\otimes y)\bullet z-x\bullet y\bullet z$$ 
Un simple raisonnement par r\'ecurrence nous conduit donc au
th\'eor\`eme suivant~: 
\th{I.1}
Pour tout espace vectoriel $V$ de dimension finie sur un corps $k$ de
caract\'eristique z\'ero, le quotient $\A$ de la big\`ebre
bitensorielle $T(T(V))$ par le bi-id\'eal ${\cal I}$ engendr\'e par la
relation~: $c=1$ admet un (unique) antipode $S$, d'ordre infini,
donn\'e par la formule~: 
$$S(x_1\otimes\cdots\otimes
x_k)=\sum_{J\subset\{1,\ldots,k-1\}}(-1)^{k+|J|}\varphi_J(\uple xk)$$ 
avec~:
$$\varphi_J(\uple xk)=x_1*_J\cdots*_Jx_k$$
o\`u $x_j*_Jx_{j+1}$ est \'egal \`a $x_j\otimes x_{j+1}$ si $j\in J$,
et \`a $x_j\bullet x_{j+1}$ si $j\notin J$, la loi $\otimes$ \'etant
suppos\'ee prioritaire devant la loi $\bullet$. 
\dem
On a v\'erifi\'e le th\'eor\`eme pour $k\le 3$. Une simple
r\'ecurrence \`a l'aide de la partie sup\'erieure du diagramme fournit
la formule donn\'ee dans l'\'enonc\'e, et on v\'erifie facilement que
la partie inf\'erieure du diagramme commute aussi. On obtient ainsi
l'expression de l'antipode sur $T(V)$, et donc sur tout $\A$ puisque
$S$ est un antimorphisme d'alg\`ebres. Reste \`a v\'erifier que le
diagramme commute pour tout \'el\'ement de $\A$. Ceci r\'esulte
directement de la proposition suivante~: [K, lemma III.3.6] 
\prop{I.2}
Soit $(B,m,\Delta,u,\varepsilon)$ une big\`ebre, $G$ une partie de $B$
engendrant $B$ en tant qu'alg\`ebre, et $S$ un anti-morphisme
d'alg\`ebres tel que le diagramme (*) commute pour tout \'element de
$G$. Alors $S$ est un antipode. 
\ndem
{\sl Fin de la d\'emonstration du th\'eor\`eme\/}~: On remarque que
pour tout $x,y\in V$ on a~: 
$$S^2(x\otimes y)=x\otimes y-(x\bullet y-y\bullet x)$$
d'o\`u on d\'eduit imm\'ediatement~:
$$S^{2p}(x\otimes y)=x\otimes y-p(x\bullet y-y\bullet x)$$
ce qui montre que l'antipode est d'ordre infini.
\qed
\alinea{I.3. L'op\'erateur de c\'esure}
On supposera dor\'enavant que le corps de base est de
caract\'eristique nulle. On d\'esigne par $S_0$ l'unique
anti-automorphisme de l'alg\`ebre $\A$ qui v\'erifie~: 
$$S_0\restr{T(V)}=-I+2u\varepsilon$$
Le compos\'e $SS_0$ est alors l'unique automorphisme de l'alg\`ebre
$\A$ qui s'exprime sur $T(V)$ par~: 
$$SS_0(x_1\otimes\cdots\otimes
x_k)=\sum_{J\subset\{1,\ldots,k-1\}}(-1)^{k+|J|+1}\varphi_J(\uple
xk)$$ 
\th{I.3}
Soit $U$ l'unique d\'erivation de l'alg\`ebre $\A$ telle que~:
$$U\restr{T(V)}=m\Delta-2I+u\varepsilon$$
Alors on a~:
$$SS_0=\exp -U$$
\dem
$U(1)=0$, et $U(x)=0$ pour tout $x\in V$. Sur les \'el\'ements de
$T(V)$ de degr\'e sup\'erieur on a~: 
$$U(x_1\otimes\cdots\otimes
x_k)=\sum_{j=1}^{k-1}(x_1\otimes\cdots\otimes
x_j)\bullet(x_{j+1}\otimes\cdots\otimes x_k)$$ 
On en d\'eduit facilement que pour $r<k$ on a~:
$$U^r(x_1\otimes\cdots\otimes x_k)=r!\sum_{1\le j_1<\cdots <j_r\le
k-1}x_1\otimes\cdots\otimes x_{j_1}\bullet
x_{j_1+1}\otimes\cdots\otimes x_{j_2}\bullet\cdots\bullet
x_{j_r+1}\otimes\cdots\otimes x_k$$ 
alors que $U^r(x_1\otimes\cdots\otimes x_k)=0$ pour $r\ge k$. Le
th\'eor\`eme s'en d\'eduit aussit\^ot. 
\qed
\cor{I.4}
Soit $C$ une sous-cog\`ebre de $T(V)$. Alors la big\`ebre ${\cal
C}=T(C)/{\cal I}\cap T(C)$ est une sous-alg\`ebre de Hopf de $\A$. 
\dem
Il est clair que $U$ respecte ${\cal C}$, donc $SS_0=\exp -U$
aussi. Comme d'autre part il est \'evident que $S_0$ a la m\^eme
propri\'et\'e, on en d\'eduit que ${\cal C}$ est stable par $S=(\exp
-U)S_0$ 
\qed
{\bf D\'efinition}~: On appelle {\sl op\'erateur de c\'esure de
l'alg\`ebre de Hopf $\A$} la d\'erivation $U$ d\'efinie dans ce
paragraphe. Cet op\'erateur associe \`a un mot la somme de toutes les
phrases obtenues en coupant ce mot en deux mots effectifs, c'est \`a
dire de longueur non nulle. 
\medskip
{\sl Remarque\/}~: une autre m\'ethode de construction de l'antipode
fait appel au produit de convolution~: on consid\`ere l'\'egalit\'e
formelle~: 
$$\eqalign{S=I^{*-1}
&=\bigl(u\varepsilon-(u\varepsilon-I)\bigr)^{*-1}\cr 
			&=\sum_{k\ge 0}(u\varepsilon-I)^{*k}}$$
On montre ensuite facilement que pour tout $v\in T(V)$ et pour tout
$k\ge 2$ on a~: 
$$(u\varepsilon-I)^{*k}(v)=(-1)^k{U^{k-1}\over (k-1)!}(v)$$
ce qui donne un sens \`a l'\'egalit\'e ci-dessus, et ce qui montre
aussi que pour tout $v\in T(V)$ on a~: 
$$S(v)=2u\varepsilon(v)-\exp {-U}(v)$$
ce qui nous permet de retrouver l'expression de l'antipode donn\'ee
par le th\'eor\`eme I.3. 
\alinea{I.4. Une famille d'\'el\'ements primitifs}
On suppose toujours que le corps est de caract\'eristique z\'ero, et
on conserve les notations du paragraphe pr\'ec\'edent. Un {\sl
\'el\'ement primitif\/} d'une alg\`ebre de Hopf est un \'el\'ement $v$
qui v\'erifie~: 
$$\Delta(v)=1\tilde\otimes v+v\tilde\otimes 1$$	
Un \'el\'ement primitif v\'erifie toujours~: $S(v)=-v$.
\smallskip\qquad
Soit $\A$ l'alg\`ebre de Hopf bitensorielle construite sur un espace
vectoriel $V$ de dimension finie. Alors tout \'el\'ement de $V$ est un
\'el\'ement primitif de $\A$. Si $x$ et $y$ appartiennent \`a $V$ on
v\'erifie ais\'ement que l'\'el\'ement~: 
$$v=x\otimes y+y\otimes x-{1\over 2}(x\bullet y+y\bullet x)$$
est primitif. Nous allons d\'ecrire un proc\'ed\'e qui permet
d'associer un \'el\'ement primitif \`a tout tenseur sym\'etrique~: 
\smallskip\qquad 
On d\'esigne par $\varphi$ la s\'erie enti\`ere d\'efinie par~:
$$\varphi(z)={1-e^{-z}\over z}=\sum_{p\ge 0}{(-1)^p\over(p+1)!}\,z^p$$
\th{I.5}
Pour tout tenseur sym\'etrique $v_0\in T(V)$ l'\'el\'ement~:
$$v=\varphi(U)(v_0)$$
est primitif.
\dem
On consid\`ere l'espace $\A[[t]]$ des s\'eries formelles \`a
coefficients dans $\A$. Cet espace est muni d'une structure
d'alg\`ebre de Hopf topologique sur $k[[t]]$ [Dr, E-K]. Pour tout
$x\in V$ on consid\`ere l'\'element~: 
$$g_x=(1-tx)^{\otimes -1}=1+\sum_{k\ge 1}t^kx^{\otimes k}$$
Un calcul imm\'ediat montre que cet \'el\'ement est {\sl de type
groupe\/}, c'est \`a dire qu'il v\'erifie~: 
$$\Delta g_x=g_x\tilde\otimes g_x$$
On en d\'eduit que~:
$$\log g_x=\sum_{k\ge 1}t^kx^{\otimes k}-{1\over 2}\bigl(\sum_{k\ge
1}t^kx^{\otimes k}\bigr)^{\bullet 2}+{1\over 3}\bigl(\sum_{k\ge
1}t^kx^{\otimes k}\bigr)^{\bullet 3}-\cdots$$ 
est primitif dans $\A[[t]]$, c'est \`a dire que pour tout entier $n\ge
1$ le terme en $t^n$ dans la somme ci-dessus est primitif dans
$\A$. Or ce terme est pr\'ecis\'ement $\varphi(U)(x^{\otimes
n})$. Pour passer du tenseur sym\'etrique $x^{\otimes n}$ \`a un
tenseur sym\'etrique quelconque on utilise le principe
d'inclusion-exclusion~: 
\smallskip\qquad
On consid\`ere dans l'ensemble $\{1,\ldots,n\}^k$ le sous-ensemble $B$
des k-uplets o\`u toutes les lettres apparaissent, et pour tout
$j\in\{1,\ldots ,n\}$ le sous-ensemble $A_j$ des k-uplets o\`u la
lettre $j$ n'appara\^\i t pas. On a alors~: 
$$\eqalign{	&\!\!\!\!\!\!\!\!(x_1+\cdots+x_n)^{\otimes k}
=\!\!\sum_{(\uple jk)\in B}\!\!x_{j_1}\otimes\cdots\otimes
x_{j_k}+\!\!\sum_{p=1}^n(-1)^{p+1}\sum_{(\uple jk)\in
A_{q_1}\!\!\cap\cdots\cap A_{q_p}}x_{j_1}\otimes\cdots\otimes
x_{j_k}\cr &=\sum_{(\uple jk)\in B}x_{j_1}\otimes\cdots\otimes
x_{j_k}+\sum_{p=1}^n(-1)^{p+1}\sum_{1\le q_1<\cdots<q_p\le n}
(x_1+\cdots\widehat{x_{q_1}}\cdots\widehat{x_{q_p}}\cdots+x_n)^{\otimes
k}\cr}$$   
Il suffit alors de faire $k=n$ et d'appliquer $\varphi(U)$ aux deux
membres de l'\'egalit\'e ci-dessus pour obtenir le th\'eor\`eme. 
\qed
\alinea{I.5. Le foncteur $\A$}
\qquad L'op\'eration qui consiste \`a associer \`a un espace vectoriel
$V$ de dimension finie son alg\`ebre de Hopf bitensorielle point\'ee
$\A_V$ est un foncteur de la cat\'egorie des espaces vectoriels de
dimension finie dans la cat\'egorie des alg\`ebres de Hopf point\'ees
(c'est \`a dire contenant un unique \'el\'ement de type groupe). En
effet si $f:V\rightarrow W$ est une application lin\'eaire entre deux
espaces vectoriels de dimension finie, $T(T(f))$ envoie l'id\'eal
${\cal I}$ de $V$ dans l'id\'eal ${\cal I}$ de $W$ et d\'efinit donc
par passage au quotient un morphisme de big\`ebres~: 
$$\A_f:\A_V\rightarrow \A_W$$
qui est donc un morphisme d'alg\`ebres de Hopf [Sw]
\smallskip\qquad
Enfin l'alg\`ebre de Hopf point\'ee $\A_V$ v\'erifie la propri\'et\'e
universelle suivante~: pour toute alg\`ebre de Hopf point\'ee $\H$ et
pour tout morphisme de cog\`ebres $\psi$ de $T(V)$ dans $\H$ il existe
un unique morphisme d'alg\`ebres de Hopf $\overline\psi$ de $\A_V$
dans $\H$ tel que le diagramme ci-dessous commute~: 
\medskip
$$\matrix{T(V)		&\inj 8 	&\A_V\cr
	\Big\downarrow\psi	&\ \ \swarrow\overline \psi &\cr
	\H		&&\cr}$$
\alinea{I.6. L'alg\`ebre de Hopf bitensorielle universelle}
\qquad On peut modifier un peu la construction ci-dessus de mani\`ere
\`a obtenir une alg\`ebre de Hopf qui v\'erifie la propri\'et\'e
universelle du \S I.5 pour toute alg\`ebre de Hopf $\H$ point\'ee ou
non~: on consid\`ere dans $\A_0[t]$ l'id\'eal $\J$ engendr\'e par
$\{c\bullet v-v\bullet c, \ v\in \A_0\}$ et par $ct-1$, et on pose~: 
$$\tilde\A=\A_0[t]/\J$$
On v\'erifie facilement que $\tilde \A$ est une alg\`ebre de Hopf, et
qu'on a une formule explicite pour l'antipode~: si $U$ d\'esigne
l'op\'erateur de c\'esure d\'efini comme pour $\A$ on a~: 
$$SS_0=C\inver\exp -U C\inver$$
o\`u $C$ d\'esigne l'automorphisme d'alg\`ebres tel que $C(v)=c\bullet
v$ pour tout $v\in T(V)$. Pour tout tenseur sym\'etrique $v\in T(V)$
c'est alors l'\'el\'ement $C\inver\varphi(U)(v)$ qui est primitif. 
\paragraphe{II. Un couplage (d\'eg\'en\'er\'e) entre deux alg\`ebres
de Hopf} 
\qquad Soit $V$ un espace vectoriel de dimension finie, et $V^*$ son
dual. Soit $\A$ (resp. $\overline\A$) l'alg\`ebre de Hopf
bitensorielle construi\-te sur $V$ (resp. $v^*$). On d\'esigne par la
m\^eme lettre $U$ les op\'erateurs de c\'esure de $\A$ et
$\overline\A$. 
\smallskip\qquad
Dans ce paragraphe on met en \'evidence un couplage d\'eg\'en\'er\'e
non nul entre les deux alg\`ebres de Hopf $\A$ et $\overline\A$. 
\th{II.1}
Il existe un couplage de Hopf non nul $<.,.>$ entre les deux
alg\`ebres de Hopf $\A$ et $\overline\A$, tel que tout $v\in
\varphi(U).S^{(2)}(V)$ est dans le noyau $\overline\A^\perp$
(resp. tout $\alpha\in \varphi(U).S^{(2)}(V^*)$ est dans le noyau
$\A^\perp$).  
\dem
Nous aurons besoin de quelques pr\'eliminaires combinatoires. A toute
phrase de $\A$ ou de $\overline\A$ on peut associer son {\sl type\/}
$J=(\uple jr)$, o\`u $j_q$ est la longueur du $q^{\hbox{\sevenrm
i\`eme}}$ mot la constituant. On appellera aussi indiff\'eremment {\sl
type\/} la partie de $\hbox{\math N}^{*2}$ constitu\'ee des $(q,j)$,
$q\le r$ et $j\le j_q$, et on notera cette partie par la m\^eme lettre
$J$. On notera $J_p$ la $p^{\hbox{\sevenrm i\`eme}}$ ligne du type
$J$. Tout type se repr\'esente par un diagramme de Young~: ci-dessous
le type $(3,4,2)$. 
$$\matrix{\bullet&\bullet&\bullet&\cr
	\bullet&\bullet&\bullet&\bullet\cr
	\bullet&\bullet&&\cr}$$
\qquad Pour toute phrase $v$ de type $J$ dont les lettres s'\'ecrivent
$(x_q)_{q\in J}$ dans l'ordre lexicographique, et pour toute partie
$L$ de $J$ on peut faire de $L$ un type de fa\c con \'evidente (en
"justifiant \`a gauche"), et consid\'erer la phrase $v_L$ de type $L$
extraite de $v$ dont les lettres sont les $(x_q,q\in L)$, rang\'ees
par ordre lexicographique. 
\smallskip\qquad
On consid\`ere sur $\hbox{\math N}^{*2}$ deux relations d'ordre
strict~: l'ordre lexicographique strict $<$ et un ordre partiel strict
$\prec$ d\'efini par~: $(x,y)\prec (x',y')$ si et seulement si $x=x'$
et $y<y'$. 
\smallskip\qquad
Soient $J$ et $K$ deux types de m\^eme cardinal. Une {\sl bonne
bijection\/} de $J$ sur $K$ est une bijection $\varphi$ telle que~: 
$$\eqalign{\varphi(x,y)\prec \varphi(x',y')	&\Longrightarrow
(x,y)<(x',y')\cr 
		(x,y)\prec (x',y')		&\Longrightarrow
		\varphi(x,y)<\varphi(x',y')\cr}$$ 
La proposition suivante d\'ecoule imm\'ediatement de la d\'efinition~: 
\prop{II.2}
Soient $J$ et $K$ deux types de m\^eme cardinal. Alors si $\varphi$
est une bonne bijection de $J$ sur $K$, $\varphi\inver$ est une bonne
bijection de $K$ sur $J$. 
\ndem
On notera $\sigma_{JK}$ l'ensemble des bonnes bijections de $J$ sur
$K$. A toute bonne bijection $\varphi$ entre deux types $J=(\uple jr)$
et $K=(\uple ks)$ de m\^eme cardinal on associe une {\sl $J-$partition
de $K$}~: 
$$K=K_1(J,\varphi)\amalg\cdots\amalg K_r(J,\varphi)$$
avec~:
$$K_p(J,\varphi)=\varphi(J_p)$$
D'apr\`es la proposition II.2 on peut associer \`a toute $J-$
partition de $K$ une $K-$ partition de $J$ et r\'eciproquement. On
dira que ces partitions de $J$ et de $K$ sont duales l'une de
l'autre. Ci-dessous nous avons repr\'esent\'e une bonne bijection
entre deux types $J$ et $K$ de cardinal $9$, et les deux partitions
duales associ\'ees. 
$$\matrix{1&2&3&4\cr
	5&6&&\cr
	7&8&9&\cr}	\hbox to 3cm{}\matrix{7&8&9\cr
						1&&\cr
						2&5&\cr
						3&4&6\cr}$$
Il nous reste encore \`a d\'efinir le {\sl nombre d'horizontalit\'e\/}
$h_\varphi$ d'une bonne bijection $\varphi$ entre deux types $J=(\uple
jr)$ et $K=(\uple ks)$ de m\^eme cardinal~: 
$$h_\varphi=\prod_{{p=1,\ldots,r\atop
q=1,\ldots,s}}\bigl(\mop{card}(\varphi(J_p)\cap K_q)\bigr)!$$ 
\prop{II.3}
$\varphi$ et $\varphi\inver$ ont m\^eme nombre d'horizontalit\'e
\dem
En appliquant la bijection $\varphi\inver$ aux ensembles
$\varphi(J_p)\cap K_q$ on a aussi~: 
$$h_\varphi=\prod_{{p=1,\ldots,r\atop
q=1,\ldots,s}}\bigl(\mop{card}(J_p\cap \varphi\inver(K_q))\bigr)!$$ 
qui est par d\'efinition \'egal \`a $h_{\varphi\inver}$.
\qed
\smallskip\qquad
On d\'efinit alors le couplage entre $\A$ et $\overline\A$ de la fa\c
con suivante~: 
\smallskip
1) L'unit\'e de $\A$ (resp $\overline\A$) est orthogonale \`a toute
phrase de $\overline\A$ (resp. $\A$) de longueur non nulle, et
$<1_{\A},1_{\overline\A}>=1$ 
\smallskip
2) Si $v$ est une phrase de $\A$ de type $J=(\uple jr)$ et $\alpha$
une phrase de $\overline\A$ de type $K=(\uple ks)$, dont les lettres
sont not\'ees $(x_q)_{q\in J}$ (resp. $(\xi_q)_{q\in K}$) alors~: 
$$<v,\alpha>=\sum_{\varphi\in\sigma_{JK}}{1\over h_\varphi}\prod_{q\in
J}<x_q,\xi_{\varphi(q)}>$$ 
ou encore, compte tenu de la proposition II.3~:
$$<v,\alpha>=\sum_{\psi\in\sigma_{KJ}}{1\over h_\psi}\prod_{q\in
K}<x_{\psi(q)},\xi_q>$$ 
En particulier~:
a) Si $v$ et $\alpha$ sont de longueurs diff\'erentes, $<v,\alpha>=0$
\smallskip
b) pour toute famille $(\uple xn)$ d'\'el\'ements de $V$ et pour toute
famille $(\uple \xi n)$ d'\'el\'ements de $V^*$ on a~: 
$${<x_1\otimes\cdots\otimes
x_n,\xi_1\otimes\cdots\otimes\xi_n>}={1\over
n!}\prod_{i=1}^n<x_i,\xi_i>$$ 
\smallskip
c) Si $v=v_1\bullet\cdots\bullet v_r$ est une phrase de $\A$ de type
$J=(\uple jr)$ de longueur $n$, chaque $v_p$ \'etant un mot de
longueur $j_p$, on a pour tout mot $\alpha$ de longueur $n$ dans
$T(V^*)\subset \overline\A$~: 
$$<v,\alpha>={n!\over j_1!\cdots j_r!}<v_1\otimes\cdots\otimes
v_k,\alpha>$$ 
\smallskip 
d) Si $\alpha=\alpha_1\bullet\cdots\bullet \alpha_s$ est une phrase de
$\overline\A$ de type $K=(\uple ks)$ de longueur $n$, chaque
$\alpha_p$ \'etant un mot de longueur $k_p$, on a pour tout mot $v$ de
longueur $n$ dans $T(V)\subset \A$~: 
$$<v,\alpha>={n!\over k_1!\cdots
k_s!}<v,\alpha_1\otimes\cdots\otimes\alpha_s>$$ 
\smallskip\qquad
Soit une phrase $v$ de $\A$, de longueur $n$, qui s'\'ecrit comme
produit de deux sous-phrases~: $v=v_1\bullet v_2$. Le type $J$ de $v$
s'\'ecrit (avec une notation \'evidente qui se passe de
commentaires)~: $J={\displaystyle{J_1\atop J_2}}$, o\`u $J_1$
(resp. $J_2$) d\'esigne le type de $v_1$ (resp. $v_2$). 
\smallskip\qquad
Soit maintenant $\alpha$ une phrase de $\overline\A$ de longueur $n$
de type $K$. Toute bonne bijection $\varphi$ de $J$ sur $K$ induit par
restriction une bonne bijection $\varphi_1$ de $J_1$ sur une partie
$K'$ de $K$, et une bonne bijection $\varphi_2$ de $J_2$ sur le
compl\'ementaire $K"$. R\'eciproquement tout couple
$(\varphi_1,\varphi_2)$ de bonnes bijections,
$\varphi_1:J_1\rightarrow K'$ et $\varphi_2:J_2\rightarrow K"$, o\`u
$(K',K")$ forme une partition de $K$ telle qu'aucun \'el\'ement de
$K"$ ne pr\'ec\`ede par l'ordre $\prec$ un \'el\'ement de $K'$, induit
une bonne bijection $\varphi$ de $J$ sur $K$, et il est clair qu'on a~: 
$$h_\varphi=h_{\varphi_1}h_{\varphi_2}$$ 
Appelons un couple $(K',K")$ comme ci-dessus un {\sl bon d\'ecoupage
de $K$\/}. Soient $(x_q)_{q\in J}$ (resp. $(\xi_q)_{q\in K}$) les
lettres qui composent $v$ (resp. $\alpha$). On a d'apr\`es ce qui
pr\'ec\`ede~: 
$$(*)\hbox to 5mm{}<v_1\bullet v_2,\,\alpha>=\sum_{(K',K")\in
D_K}\sum_{{\varphi_1\in \sigma_{J_1K'}\atop \varphi_2\in
\sigma_{J_2K"}}}{1\over h_{\varphi_1}h_{\varphi_2}} 
\prod_{q\in J_1}<x_q,\xi_{\varphi_1(q)}>\prod_{q\in
J_2}<x_q,\xi_{\varphi_2(q)}>$$ 
o\`u $D_K$ d\'esigne l'ensemble des bons d\'ecoupages de $K$. D'autre
part, on voit facilement que~: 
$$\Delta\alpha=\sum_{(K',K")\in
D_K}\alpha_{K'}\tilde\otimes\alpha_{K"}$$ 
d'o\`u gr\^ace \`a (*)~:
$$<v_1\bullet v_2,\,\alpha>=<v_1\tilde\otimes v_2,\,\Delta\alpha>$$
La d\'emonstration de la deuxi\`eme assertion du th\'eor\`eme est
alors imm\'ediate, compte tenu du r\^ole sym\'etrique jou\'e par $\A$
et $\overline\A$, mis en \'evidence par les propositions II.2 et
II.3. Le th\'eor\`eme II.1 est donc d\'emontr\'e. 
\qed
\cor{II.4}
L'orthogonal $\overline\A^\perp$ (resp. $\A^\perp$) est un bi-id\'eal
de $\A$ (resp. $\overline\A$) 
\ndem
\prop{II.5}
Les antipodes $S$ et $\overline S$ des alg\`ebres de Hopf respectives
$\A$ et $\overline\A$ sont adjoints l'un de l'autre pour le couplage. 
\dem 
Il faut d\'emontrer que pour toute phrase $v\in \A$ et
$\alpha\in\overline\A$ on a~: 
$$<Sv,\,\alpha>=<v,\,\overline S\alpha>$$
Si $v$ et $\alpha$ sont des mots cela d\'ecoule des deux \'egalit\'es~:
$$\eqalign{<S_0v,\,\alpha>	&=<v,\,S_0\alpha>\cr
		<U(v),\,\alpha>	&=<v,\,U(\alpha)>\cr}$$
et du th\'eor\`eme I.3. Si maintenant $v_1$ et $v_2$ sont deux phrases
telles que pour tout mot $\alpha$ on ait
$<Sv_i,\,\alpha>=<v_i,\,\overline S\alpha>,\,i=1,2$, on a~: 
$$\eqalign{<S(v_1\bullet v_2),\,\alpha>	&=<Sv_2\bullet Sv_1,\,\alpha>\cr
					&=<Sv_2\tilde\otimes
					Sv_1,\,\Delta\alpha>\cr 
					&=<v_1\tilde\otimes
					v_2,\,(\overline
					S\tilde\otimes\overline
					S)\Delta^{\smop{op}}\alpha>\cr
					&=<v_1\tilde\otimes
					v_2,\,\Delta \overline
					S\alpha>\cr 
					&=<v_1\bullet
					v_2,\,S\alpha>\cr}$$ 
ce qui montre que l'\'egalit\'e a lieu pour toute phrase $v$ et tout
mot $\alpha$. R\'ep\'etant ce proc\'ed\'e d'extension du c\^ot\'e dual
on obtient la proposition. Le couplage est donc un {\sl couplage
d'alg\`ebres de Hopf\/}. 
\qed
\cor{II.6}
L'orthogonal $\overline\A^\perp$ (resp. $\A^\perp$) est un id\'eal de
Hopf de $\A$ (resp. $\overline\A$) 
\ndem
Nous allons maintenant montrer que le couplage est d\'eg\'en\'er\'e,
et mettre en \'evidence un id\'eal de Hopf explicite non nul $\J$
(resp. $\overline\J$) contenu dans l'orthogonal $\overline\A^\perp$
(resp. $\A^\perp$). 
\lemme{II.7}
pour tout $n\ge 2$, pour tout $\uple xn\in V$ (resp. $\uple\xi n\in
V^*$) et pour tout mot $\alpha\in T(V^*)\subset\overline\A$
(resp. $v\in T(V)\subset\A$) on a~: 
$$<\varphi(U)(x_1\otimes\cdots\otimes x_n),\alpha>=0
\hbox to 15mm{\hfill(resp. }<v,\,\varphi(U)(\xi_1\otimes\cdots\otimes
\xi_n)>=0\hbox{)}$$ 
\dem
Les espaces $V$ et $V^*$ jouant un r\^ole sym\'etrique, il suffit de
montrer la premi\`ere assertion du lemme. On a~: 
$$<\varphi(U)(x_1\otimes\cdots\otimes
x_n),\alpha>=\sum_{k=1}^n{(-1)^k\over
k+1}A_{n,k}<x_1\otimes\cdots\otimes x_n,\alpha>$$ 
avec~:
$$A_{n,k}=\sum_{1\le j_1<\cdots<j_k<n}{n!\over
j_1!(j_2-j_1)!\cdots(n-j_k)!}$$ 
On introduit la famille de polyn\^omes $P_n$ d\'efinie pour $n\ge 1$ par~:
$$P_n(t)=\sum_{k=0}^{n-1}t^kA_{n,k}$$
On a pour tout mot $\alpha\in T^n(V^*)\subset\overline\A$~:
$$<e^{tU}(x_1\otimes\cdots\otimes
x_n),\alpha>=P_n(t)<x_1\otimes\cdots\otimes x_n,\alpha>$$ 
et donc le lemme est vrai si et seulement si pour tout $n\ge 2$ on a~:
$$\int_{-1}^0 P_n(t)\,dt=0$$
L'\'egalit\'e~:
$${U^k\over k!}(x_1\otimes\cdots\otimes
x_n)=\sum_{j=1}^{n-1}x_1\otimes\cdots\otimes x_j\bullet {U^{k-1}\over
(k-1)!}(x_{j+1}\otimes\cdots\otimes x_n)$$ 
conduit \`a la relation de r\'ecurrence~:
$$P_n(t)=1+t\sum_{k=1}^{n-1}C_n^kP_k(t),\ \ n\ge 2$$
avec bien s\^ur $P_1=1$. On pose $P_0=0$, ce qui nous permet
d'\'ecrire pour $n\ge 1$~: 
$$(1+t)P_n(t)=1+t\sum_{k=0}^n C_n^k P_k(t)$$
En consid\'erant une ind\'etermin\'ee suppl\'ementaire $z$ on a donc~:
$$\eqalign{(1+t)\sum_{n\ge 0}{P_n(t)\over n!}z^n
&=e^z-1+t\sum_{n\ge 0}\sum_{k=0}^n P_k(t){z^n\over k!(n-k)!}\cr 
						&=e^z-1+t\sum_{k,l\ge
						0}P_k(t){z^k\over
						k!}{z^l\over l!}\cr 
						&=e^z-1+te^z\sum_{k\ge
						0}P_k(t){z^k\over
						k!}\cr}$$ 
soit encore~:
$$\sum_{n\ge 0}{P_n(t)\over n!}z^n={e^z-1\over 1-t(e^z-1)}$$
Int\'egrant les deux membres entre $-1$ et $0$ on a donc~:
$$\int_{-1}^0\sum_{n\ge 0}{P_n(t)\over n!}z^n\,dt=\bigl[-\log
(1-t(e^z-1))\bigr]_{-1}^0=z$$ 
d'o\`u $\displaystyle\int_{-1}^0 P_n(t)\,dt=0$ pour $n=0$ et $n\ge 2$,
ce qui d\'emontre le lemme 
\qed
On note $S^k(V)$ l'espace des tenseurs sym\'etriques homog\`enes de
degr\'e $k$, et on pose~: 
$$S^{(2)}(V)=\bigoplus_{k\ge 2}S^k(V)$$
D'apr\`es le \S\ I.4, l'ensemble des $\varphi(U).v,v\in S^{(2)}(V)$
est form\'e d'\'el\'ements primitifs de \A, et donc l'id\'eal $\J$
engendr\'e par cette famille est un id\'eal de Hopf. 
\cor{II.8}
L'id\'eal de Hopf $\J$ (resp. $\overline\J$) est contenu dans
$\overline\A^\perp$ (resp. $\A^\perp$) 
\dem
Il suffit d'apr\`es le corollaire II.4 de montrer que $\varphi(U)(v)$
appartient \`a $\overline\A^\perp$ pour tout $v\in S^{(2)}(V)$. Or cet
\'el\'ement est primitif (th\'eor\`eme I.5) et
$<\varphi(U)(v),\,\alpha>=0$ pour tout mot $\alpha$, d'apr\`es le
lemme II.5. Soient maintenant deux phrases $\alpha_1$ et $\alpha_2$ de
$\overline\A$ telles que~: 
$$<\varphi(U)(v),\,\alpha_i>=0,\ i=1,2$$
Alors~:
$$\eqalign{<\varphi(U)(v),\,\alpha_1\bullet\alpha_2>
&=<\Delta\varphi(U)(v),\, \alpha_1\tilde\otimes\alpha_2>\cr 
&=<\varphi(U)(v),\,\alpha_1><1,\,\alpha_2>+<1,\,\alpha_1>
<\varphi(U)(v),\alpha_2>\cr   
&=0\cr}$$
d'o\`u on d\'eduit que $<\varphi(U)(v),\alpha>=0$ pour toute phrase
$\alpha$. Ceci d\'emontre le th\'eor\`eme II.1. 
\qed
\paragraphe{R\'ef\'erences}
\bib{Ab}E. Abe, {\sl Hopf algebras\/}, Cambridge (1980)
\bib{Dr}V.G. Drinfeld, {\sl Quantum Groups\/}, Proceedings ICM
(Berkeley 1986) {\bf 1}, 798-820. 
\bib{E-K}P. Etingof \& D. Kazhdan, {\sl quantization of Lie
bialgebras, I\/}, pr\'epubl. univ. Harvard (q-alg 9506005) 
\bib{K}Ch. Kassel, {\sl Quantum Groups}, Springer (1995)
\bib{R}N. Reshetikhin, {\sl Quantization of Lie bialgebras\/},
Int. Math. res. Notices {\bf 7} (1992) 
\bib{Sw}M.E. Sweedler, {\sl Hopf algebras}, Benjamin (1969) 

\bye